\documentclass{elsart}

\usepackage{amssymb}
\usepackage{amsmath, amsfonts, mathrsfs}
\usepackage{indentfirst}
\usepackage{latexsym,fancyhdr,fancybox}

\begin{document}
\begin{frontmatter}



\title{Difference between Devaney chaos associated with two systems}


\author[Changchun]{Bingzhe Hou\thanksref{NSFC}}
\ead{abellengend@163.com}
\author[Shanghai,Nanjing]{Xianfeng Ma\corauthref{cor}\thanksref{PSRP}}
\corauth[cor]{Corresponding author.} \ead{xianfengma@gmail.com}
\author[Changchun]{Gongfu Liao\thanksref{NSFC}}
\ead{liaogf@email.jlu.edu.cn}
\address[Changchun]{Institute of Mathematics, Jilin University, Changchun 130012, PR China}
\address[Shanghai]{Department of Mathematics, East China University of Science and
Technology, Shanghai 200237, China}
\address[Nanjing]{School of Mathematics and Computer Science, Nanjing Normal
University, Nanjing 210097, PR China}

\thanks[PSRP]{Supported by a grant from Postdoctoral Science Research Program of Jiangsu Province (No.
0701049C) and Specialized Research Fund for Outstanding Young
Teachers of East China University of Science and Technology.}
\thanks[NSFC]{Supported by NSFC (No. 10771084).}
\begin{abstract}

We discuss the relation between Devaney chaos in the base system and
Devaney chaos in its induced hyperspace system. We show that the
latter need not imply the former. We also argue that this
implication is not true even in the strengthened condition.
Additionally we give an equivalent condition for the periodically
density in the hyperspace system.
\end{abstract}


\end{frontmatter}

\section{Introduction}
Let $(X,d)$ be a compact metric space with metric $d$ and $f:
X\rightarrow X$ be continuous. Then $(X,f)$ is called a compact
system.
For every positive integer $n$, we define $f^n$ inductively
by $f^n=f\circ f^{n-1}$, where $f^0$ is the identity map on $X$.

A compact system $(X,f)$ is called  Devaney chaos \cite{Devaney} if
it satisfies the following three conditions:\\
\indent (1) $f$ is transitive, i.e., for every pair $U$ and $V$ of
non-empty open subsets of $X$ there is a non-negative integer $k$
such that $f^{n}(U)\cap V\neq\emptyset$;\\
\indent (2) $f$ is periodically dense, i.e., the set of periodic
points of $f$ is dense in $X$;\\
\indent (3) $f$ is sensitive, i.e., there is a $\delta > 0$ such
that, for any $x\in X$ and any neighborhood $V$ of $x$, there is a
non-negative integer $n$ such that $d(f^n(x),f^n(y))>\delta$.

It is worth noting that the conditions (1) and (2) imply that $f$ is
sensitive if $X$ is infinite \cite{Silverman, Banks}. This means the
condition (3) is redundant in the above definition.

In the works \cite{Roman-Flores} Rom\'{a}n-Flores investigated a
certain hyperspace system $(\mathcal{K}(X),\overline{f})$ associated
to the base system $(X,f)$, where
$\overline{f}:\mathcal{K}(X)\rightarrow\mathcal{K}(X)$ is the nature
extension of $f$ and $\mathcal{K}(X)$ is the family of all non-empty
compact sets of $(X,d)$ endowed with the Hausdorff metric induced by
$d$. He presented a fundamental question: Does the chaoticity  of
$(X,f)$ (individual chaos) imply that of
$(\mathcal{K}(X),\overline{f})$ (collective chaos)? and conversely?

As a partial response to this question, Rom\'{a}n-Flores
\cite{Roman-Flores} discussed the transitivity of the two systems,
and showed that the transitivity of $\overline{f}$ implies that of
$f$, but the converse is not true. Fedeli \cite{Fedeli} showed that
the periodically density of $f$ implies that of $\overline{f}$.
Banks  gave an example which has a dense set of periodic points in
the hyperspace system but has none in the base system \cite{Banks}.
Gongfu Liao showed that there is an example on the interval which is
Devaney chaos in the base system while is not in its induced
hyperspace system \cite{liaogf}.

In this paper, we show that $\overline{f}$ being Devaney chaos need
not imply $f$ being Devaney chaos. Further, $f$ need not be Devaney
chaos even if $\overline{f}$ is mixing and periodically dense. This
answers the question posed by Rom\'{a}n-Flores and Banks
\cite{Banks}.

\section{Preliminaries}

$(X,f)$ is a compact system. If $Y\subset X$ and $Y$ is
$f$-invariant, i.e. $f(Y)\subset Y$, then $(Y, f|_{Y})$ is called
the subsystem of $(X,f)$ or $f$, where $f|_{Y}$ is the restriction
of $f$ on $Y$.

A subset $A\subset X$ is minimal if  the closure of the orbit of any
point $x$ of $A$ is $A$, i.e. $\overline{Orb(x)}= A$, for all $x\in
A$. An equivalent notion is that $A$ has no proper $f$-invariant
closed subset.

A point $x\in X$ is said to be an almost periodic point if for any
$\varepsilon>0$ exists $N\in \mathbb{N}$ such that for any integer
$q\geq 1$  exists an integer $r$ with property that $ q\leq r <N+q$
and $f^{r}(x)\in B(x,\varepsilon)$, where $B(x,\varepsilon)=\{y\in
X: d(x,y)<\epsilon\}$. A point $x\in X$ is an almost periodic point
if and only if $\overline{Orb(x)}$ is minimal \cite{Zuoling Zhou}.

$f$ is weakly mixing if the Cartesian product $f\times f$ is
transitive, that is for any non-empty open set $U_1, U_2,V_1$ and
$V_2$, there is a positive integer $n$ such that $f^n(U_i)\cap
V_i\neq\emptyset$, $i=1,2$. $f$ is  mixing if for any non-empty open
sets $U$ and $V$, there is a positive integer $N$ such that for all
$n\geq N$, $f^n(U )\cap V \neq\emptyset$

Let $S=\{0,1,\cdots ,k-1\}$, $\Sigma_{k}=\{x=x_{0}x_{1}\cdots:
 x_{i}\in S,i=0,1,\cdots\}$.
 Define $\rho:\Sigma_{k}\times\Sigma_{k}\rightarrow R$ as follows: for any $x$, $y \in
 \Sigma_{k}$,  $x=x_{0}x_{1}\cdots$, $y=y_{0}y_{1}\cdots$,
 \begin{align*}
 \rho(x,y)=
 \begin{cases}
 \quad 0,&\text{if} ~x=y,\\
 \dfrac{1}{m+1}, &\text{if} ~x \not= y.
 \end{cases}
 \end{align*}
 where $ m=\min\{n\geq1|x_{n}\neq y_{n}\}$.
$\rho$ is a metric over $\Sigma_{k}$, and $(\Sigma_{k},\rho)$ is a
compact metric space and is called a symbol  space.

Define $\sigma:\Sigma_{k}\rightarrow\Sigma_{k}$ as
$\sigma(x)=x_{1}x_{2}\cdots$, where
$x=x_{0}x_{1}\cdots\in\Sigma_{k}$. $\sigma$ is continuous and is
called the shift map on $\Sigma_{k}$.

 $A$ is called a block over $S$ if it is a finite permutation of
 some elements of $S$. For $A=a_{0} \cdots a_{n-1}$, where $a_{i}\in
S$, $0\leq i \leq n-1$,  the length of $A$ is defined as $n$, i.e.
$|A|=n$. Denote  $[A]=\{x\in
\Sigma_{k}:~x_{i}=a_{i},~i=0,1,\cdots,n\},$ $[A]$ is called the
cylinder set over   $\Sigma_{k}$ , it is an open and closed set.
Let
 $B=b_{0} \cdots b_{m-1}$ be another block, then the concatenation $C=AB=a_{0} \cdots a_{n-1} b_{0}
\cdots b_{m-1}$ is a new block. $B$ is said to be admissible in $A$
if there exists $i\geq 0$  such that
 $b_{j}=a_{i+j},j=0,1,\cdots,m-1 .$  For $c\in S$,
denote $c^{[n]}$ as a $n$ length permutation of $c$ (For example,
$2^{[3]}=222$ ), and $c^{[\infty]}=ccc\cdots$ as an infinite
permutation of $c$.

\begin{lem}
        $ b_{0}b_{1}\cdots
          b_{n}\cdots$ is an almost periodic point of a symbol space if and only if
        for any $j>0$ exists $N>0$ such that for any $i>0$ symbol block $ b_{0}b_{1}\cdots
       b_{j}$ appears in $ b_{i}b_{i+1}\cdots b_{i+N}$ .
\end{lem}
\begin{lem}
\label{fhqhh}
 Let $(Y, \sigma)$ be a subsystem of  a symbol system $(\Sigma_{k}, \sigma)$,
 then $(Y, \sigma)$ is mixing
 if and only if
 for any two admissible blocks $b_{0}b_{1}\cdots b_{s}$ and $c_{0}c_{1}\cdots c_{t}$
 in $Y$,
 there exists $N>0$ such that for any $n>N$ exists $n-1$ length admissible block $w_{1}w_{2}\cdots
 w_{n-1}$ in $Y$ with the property that $b_{0}b_{1}\cdots b_{s}w_{1}w_{2}\cdots w_{n-1}c_{0}c_{1}\cdots
 c_{t}$ is an admissible block of $Y$.
\end{lem}

Let $(X,\mathcal{J})$ be a topological space and $\mathcal{K}(X)$ be
the set of all the non-empty compact sets of $X$. $\mathcal{K}(X)$
is called the hyperspace of $(X,\mathcal{J})$ when endowed with the
Vietoris topology $\mathcal{J}_{\mathscr{V}}$ whose base consists of
sets of the form
 \begin{equation*} \mathscr{B}(G_{1},G_{2},\cdots
,G_{n}) =\{K\in\mathcal{K}(X):K\subset
\overset{n}{\underset{i=1}{\cup}}G_{i}~\text{and}~ K\cap G_{i}\neq
\emptyset,\,\,1\leq i\leq n \},
\end{equation*}
where $G_{1},G_{2},\cdots ,G_{n}$ are non-empty open subsets of $X$.

$(X,d)$ is a metric space, $A\in \mathcal{K}(X)$,
``$\epsilon$-dilatation of $A$" is defined as the set
$N(A,\epsilon)=\{x\in X:d(x,A)<\epsilon\}$, where
$d(x,A)=\inf\limits_{a\in A}d(x,a)$.
The Hausdorff separation
$\rho(A,B)$ is $\rho(A,B)=\inf\{\epsilon>0:A\subset
N(B,\epsilon)\}$.
The Hausdorff mertic on $\mathcal{K}(X)$  is
$H_{d}(A,B)=\max\{\rho(A,B),\rho(B,A)\}$.
 It is well known   that the topology induced by the
 Hausdorff metric $H_{d}$ on $\mathcal{K}(X)$ coincides with the
 Vietoris topology $\mathcal{J}_{\mathscr{V}}$ \cite{Klein}.

The extension $\overline{f}$ of $f$ to $\mathcal{K}(X)$ is defined
as $\overline{f}(A)=\{f(a):a\in A\}$. It is easy to show that $f$
being continuous implies $\overline{f}$ being continuous. If $(X,
f)$ is compact system, then $(\mathcal{K}(X), \overline{f})$ is also
a compact system \cite{Klein}. $(\mathcal{K}(X), \overline{f})$ is
said to be the hyperspace system induced by the base system $(X,
f)$.

The following two theorems can be found in
\cite{Banks,liaogf,Peris,gurongbao}
\begin{thm}
\label{chrdj} Let $(\mathcal{K}(X),\overline{f})$ be the hyperspace
system induced by the compact system $(X,f)$, then the following are
equivalent.\\
\indent (i) $f$ is weakly mixing,\\
\indent(ii) $\overline{f}$ is weakly mixing,\\
\indent(iii) $\overline{f}$ is transitive.
\end{thm}
\begin{thm}
\label{mixing}
  Let $(\mathcal{K}(X),\overline{f})$ be the hyperspace
system induced by the compact system $(X,f)$, then $f$ being mixing
is equivalent to $\overline{f}$ being mixing.
\end{thm}

\section{Periodically dense}

\begin{thm}
\label{chbsh} If $(\mathcal{K}(X), \overline{f})$ is the hyperspace
system induced by the compact system $(X,d)$, then $\mathcal{K}(X)$
has a dense set of periodic points if and  only if for any non-empty
open subset $U$ of $X$ there exists a compact subset $K\subset U$
and an integer $n>0$ such that $f^{n}(K)=K$.
\end{thm}

\noindent
 \textbf{Proof}

 $\Longrightarrow$ Since $\mathcal{K}(X)$ has a dense set of periodic
 points, every non-empty open subset of $\mathcal{K}(X)$ has at least one periodic point.
Further, every basic element has at least one periodic point. Let
$U$ be an arbitrary open subset of $X$. Then by the definition of
the Vietoris topology $\mathscr{B}(U)$ is a basic element of
$\mathcal{K}(X)$.
 According to the condition, there exists $K\in\mathcal{K}(X)$ and an integer
$n>0$  such that $\overline{f}^{n}(K)=K$. Since
$\overline{f}^{n}=\overline{f^{n}}$, then $\overline{f^{n}}(K)=K$.
Hence, there is a $K\subset U$ in $X$ such that $f^{n}(K)=K$.

$\Longleftarrow$ Let $U_{1},U_{2},\ldots,U_{n}$ be the open subsets
of $X$, then
$$\mathscr{B}(U_{1},U_{2},\ldots,U_{n})=
\{W\in\mathcal{K}(X):W\subset\overset{n}{\underset{i=1}\cup}{U_{i}}~\text{
and}~W\cap U_{i}\neq \emptyset,  \,\, 1\leq i\leq n\}$$ is a basic
element   of $\mathcal{K}(X)$.
 According to the condition, there
exists $K_{i}\subset U_{i}, m_{i}>0$ such that
$f^{m_{i}}(K_{i})=K_{i},  \forall \,\, 1\leq i\leq n ,$  Let $m$ be
the least common multiple of $\{m_{i}\}$ and
$K=\overset{n}{\underset{i=1}\cup}{K_{i}}$.
 Then $f^{m}(K)=K$, and
 $K\subset\overset{n}{\underset{i=1}\cup}{U_{i}},\quad K\cap
U_{i}\neq \emptyset,  \,\, 1\leq i\leq n$. Therefore,
$K\in\mathscr{B}(U_{1},U_{2},\ldots,U_{n})$.

This means that every element of $\mathcal{K}(X)$ has a
$m$-invariant set of $X$ in it. Further, every open set of
$\mathcal{K}(X)$ has a periodic point in it.
 Then  $\mathcal{K}(X)$ has a dense set of periodic points.

\begin{rem}
         If $X$ is  a topology space,
         and $\mathcal{P}(X)$ is the family of all non-empty subset
         of it, then with the similar argument, it is not hard to see
         the above result is also hold for the Vietoris
         topology.
\end{rem}

\section{Devaney chaos}

If $a$ is an almost periodic point of $\Sigma_{2}$ and
$X=\overline{orb(a)}$, then $(X, \sigma)$ is usually a mixing
non-trivial minimal subsystem of $(\Sigma_{2}, \sigma)$, which
depends on the selection of the point $a$. There are many such
examples \cite{Qinjie}. From these examples we take an arbitrary one
as our research object and this do not affect the proof and our
results.

We construct a subsystem $(\widetilde{X}, \sigma)$ of $(\Sigma_{3},
\sigma)$ associated to the selected $(X, \sigma)$ as follows.

For $b=b_{0}b_{1}\cdots b_{n}\cdots\in X$, denote
$$
X_{b} = \{ x = x_{0}x_{1}\cdots\in \Sigma_{3}: \, \exists
n_{i}\rightarrow\infty\,\text{such that}\,x_{n_{i}}= b_{i},
\,\text{and}\,x_{n} = 2 \,\text{if}\,n\not= n_{i}\},
$$
then $X_{b}$ is called the extension of $b$ in $\Sigma_{3}$.
Denote
 $\widetilde{X} = \overline{\bigcup_{b\in X}X_{b}}.$
It is easy to see that $\widetilde{X}$ is $\sigma$-invariant set of
$\Sigma_{3}$. Therefore, $(\widetilde{X},\sigma)$ is a subsystem of
$(\Sigma_{3},\sigma)$.

\begin{prop}
\label{xqhh}
         $(\widetilde{X},\sigma)$  is mixing.
\end{prop}

\noindent
 \textbf{Proof}
 Let  $U=[x_{0}x_{1}\cdots x_{n_{1}}]$ and $V=[y_{0}y_{1}\cdots
y_{n_{2}}]$ be the cylinder sets of $\widetilde{X}$. They are two
basic elements of $\widetilde{X}$ and both are closed sets.

By the construction of $\widetilde{X}$, there exists two integer
sequences $\{i_{l}\},\{j_{m}\},i_{l},j_{m}\geq-1,l,m\geq-1$~ such
that $x_{i_{l}}=b_{l},y_{j_{m}}=c_{m}$, and
 $U_{1}=[b_{0}b_{1}\cdots b_{l}\cdots b_{s}],
V_{1}=[c_{0}c_{1}\cdots c_{m}\cdots b_{t}]$, $s\leq n_{1}, t\leq
n_{2}$,
 are cylinder sets of $X$
 (We assign $U=[2^{n_{1}}]$ for $i_{l}=-1$, and $V=[2^{n_{2}}]$ for $j_{m}=-1$).
  So
 $u=b_{0}b_{1}\cdots b_{l}\cdots b_{s}, \quad v=c_{0}c_{1}\cdots
c_{m}\cdots b_{t}$  are the admissible blocks of $X$.

For $i_{l}\geq 0$, $j_{m}\geq 0$, $\forall \,\, l,m$, from Lemma
 \ref{fhqhh}, there exists $N\in \mathbb{N}$ such that for any $n\geq
 N$ there exists an admissible block $w=d_{1}d_{2}\cdots d_{n-1}$ of $X$
 with the property that $uwv$ is an admissible block of $X$.
Since $U,V$ can be written as
 \begin{align*}
&U=[x_{0}x_{1}\cdots x_{i_{s}}2^{[r]}],\,\, r=n_{1}-i_{s}   \\
&V=[2^{[j_{0}]}y_{j_{0}}y_{j_{0}+1}\cdots y_{n_{2}}],
\end{align*}
according to the construction of $\widetilde{X}$, obviously,
$$x_{0}x_{1}\cdots x_{i_{s}}2^{[r]}w2^{[j_{0}]}y_{j_{0}}y_{j_{0}+1}\cdots y_{n_{2}}$$
is an admissible block of $\widetilde{X}$.
Then
$$[x_{0}x_{1}\cdots x_{i_{s}}2^{[r]}w2^{[j_{0}]}y_{j_{0}}y_{j_{0}+1}\cdots y_{n_{2}}]
\subset \widetilde{X}$$ is the basic element of $\widetilde{X}$.
This means that $\sigma^{n}(U)\cap V\neq \emptyset,\quad \forall \,
n\geq N.$

For $i_{l}$ or $j_{m}$ is $-1$, we just need to add enough symbols 2
in the admissible blocks  $U$ and $V$.

 Therefore,
$(\widetilde{X},\sigma)$ is mixing.

\begin{cor}
\label{mmxing} $(\mathcal{K}(\widetilde{X}),\overline{\sigma})$ is
mixing.
\end{cor}
This is a direct result of Proposition \ref{xqhh} and Theorem
\ref{mixing}.

\begin{prop}
\label{wyzhqd}
       $(\widetilde{X},\sigma)$ is not minimal and has only one periodic point.
       The periodic point is the unique fixed point of $\widetilde{X}$.
\end{prop}

 \textbf{Proof}
 From the construction,
 it is easy to see that $X\subset \widetilde{X}$ and $X$ is minimal proper set of
 $\Sigma_{3}$.
 By the definition of minimal, $\widetilde{X}$ is not minimal.
 We now prove $\widetilde{X}$ has only one fixed point $2^{[\infty]}$.
 It is the unique periodic point of $\widetilde{X}$.

Suppose that there exists $x\in \widetilde{X}$ and $n>0$ such that
$\sigma^{n}(x)=x$, then $x$ can be written as  an infinite
permutation of block $x_{0}x_{0}\cdots x_{n-1}$,
$$x=x_{0}x_{1}\cdots x_{n-1}\cdots x_{0}x_{1}\cdots x_{n-1}\cdots.$$
If $x_{i}=2$ for all $i$, $0\leq i\leq n-1$, then $x=2^{[\infty]}$.
Obviously, $x$ is a fixed point of $\widetilde{X}$.

If there exists a finite integer sequence $\{i_{l}\}$, $0\leq
i_{l}\leq n-1$, $0\leq l\leq m$, $0\leq m\leq n-1$ such that
$x_{i_{l}}\neq 2$, then $x_{i_{0}}x_{i_{1}}\cdots x_{i_{m}}$ is an
admissible block of $X$, we can write it for convenience  as
$a_{0}a_{1}\cdots a_{m}$.
 Accordingly,
  $a=a_{0}a_{1}\cdots
a_{m}\cdots a_{0}a_{1}\cdots a_{m}\cdots$
 is an infinite permutation
of $a_{0}a_{1}\cdots a_{m}$,
then $a$ is a $(m+1)$-periodic point.

Since $X$ is a non-trivial minimal set, so $a\not\in X$, and then
$x\not\in X_{a}$. Further, it is not hard to see that $x\not\in
\widetilde{X}$. This is a contradiction.

\begin{rem}
\label{dense}
        The result implies that $\widetilde{X}$ is not periodically dense.
\end{rem}

\begin{prop}
\label{chkjchm}
 The hyperspace system $(\mathcal{K}(\widetilde{X}),\overline{\sigma})$
 induced by $(\widetilde{X},\sigma)$ is periodically dense.
\end{prop}

\textbf{Proof} According to Theorem \ref{chbsh}, we just need to
prove that every open set $U$ of $X$ has a $m$-invariant compact
subset $K$, i.e. $\sigma^{m}(K) = K$. Further, it is enough to prove
that every basic element of $X$ has a $m$-invariant compact subset.

 Let $U = [a_{0}a_{1}\cdots a_{n-1}]$ be a cylinder set of $\widetilde{X}$.

If $x_{i}=2$ for all  $i$, $0\leq i \leq n-1$,  then the single set
$\{2^{[\infty]}\}$ ia a $\sigma$-invariant closed set of $U$. It is
a 1-invariant set of $\sigma$ and satisfies the condition.

Suppose that $a_{0}a_{1}\cdots a_{n-1}$ have symbols which are not 2
and the number of which is $j$, $1\leq j \leq n$. On the basis of
the definition of $\widetilde{X}$, there exists a point
$b=b_{0}b_{1}\cdots$ in $X$ satisfied: the first $j$ symbols in
order are just the symbols of $a_{0}a_{1}\cdots a_{n-1}$ which are
not 2.

Since $b$ is an almost periodic point, there exists $N > 0$ and
$n_{k}\rightarrow\infty$ such that
 $ j\leq n_{k+1} - n_{k}\leq N+j$
and the first $j$ symbols of $\sigma^{n_{k}}(b)$ and $b$ are
completely same.

Let $m = n+N$ and define $\bar{b}= {\bar b_{0}}{\bar b_{1}}\cdots\in
\Sigma_{3}$  as follows. For any $s \geq 0$,
\begin{align*}
{\bar b_{s}}=
\begin{cases}
  a_{i},          &\ \text{if}~ s = km +i,    \quad     \text{where} \ k \geq 0,\,\, 0\leq i < n,\\
 b_{n_{k}+j+i},  &\ \text{if}~ s = km +n +i,   \quad   \text{where} \ k \geq 0, \,\,0\leq i < n_{k+1}- n_{k} - j,\\
 2,              &\ \text{other}.
 \end{cases}
\end{align*}
It is not hard to see that $\bar{b}\in X_{b}\subset \widetilde{X}$.

On the other hand, obviously, for any $k \geq
 0$, $\sigma^{km}(\bar{b})\in U$.

 Let $K $ be the $\omega$-limit set of the point $\bar b$ under the action of
 $\sigma^{m}$, i.e. $K = \omega(\bar{b}, \sigma^{m})$,
 then $K\subset U$ is a closed set and $\sigma^{m}(K) = K$.

The proof of Proposition \ref{chkjchm} is completed.

\begin{thm}
\label{chchndch}  The hyperspace system being periodically dense
need not imply the base system being periodically dense.
\end{thm}
\textbf{Proof} It is a direct result of Proposition \ref{chkjchm}
and Remark \ref{dense}.

\begin{cor}
The hyperspace system being mixing and periodically dense need not
imply the base system being periodically dense.
\end{cor}

\textbf{Proof} We can get the result from Proposition \ref{xqhh},
\ref{wyzhqd}, \ref{chkjchm} and Corollary \ref{mmxing}.
\begin{thm}
The hyperspace system being Devaney chaos need not imply the base
system being Devaney chaos.
\end{thm}

\textbf{Proof} Since mixing implies weakly mixing, and weakly mixing
is equivalent to the transitivity by Theorem \ref{chrdj}, then the
hyperspace system being mixing and periodically dense implies it
being Devaney chaos. However, from Theorem \ref{chchndch}, the base
system need not be periodically dense, then it need not be Devaney
chaos.

\section{Conclusion}

In summary, we showed some diffence between the chaoticity of a
compact system and the chaoticity of its induced hyperspace system.
Devaney chaos in the hyperspace system need not implies Devaney
chaos in its base system. Even if the hyperspace system is mixing,
the former need not implies the latter. This kind of investigation
should be useful in many real problems, such as in ecological
modeling, demographic sciences associated to migration phenomena and
numerical simulation, etc.



\end{document}